\documentclass[11pt]{article}
\usepackage{verbatim}
\usepackage{amsfonts,amsmath,amsthm}

\usepackage{amssymb}
\usepackage{float}
\usepackage{epsfig}
\usepackage{amsmath}
\usepackage[english]{babel}
\usepackage{a4}
\usepackage{amstext}
\usepackage{mathrsfs}
\usepackage{amsthm}
\usepackage{fancyhdr}
\usepackage{latexsym}
\usepackage{amsmath}
\usepackage{amsfonts}
\usepackage{amssymb}
\usepackage{fancyhdr}

\usepackage[usenames]{color}
\definecolor{red}{rgb}{1.0,0.0,0.0}

\definecolor{blu}{rgb}{0.0,0.0,1.0}

\def\to{\longrightarrow}

\def\norm{{\| \kern -.05em | }}

\newtheorem{theorem}{Theorem}[section]
\newtheorem{proposition}[theorem]{Proposition}

\newtheorem{corollary}[theorem]{Corollary}
\newtheorem{assumption}[theorem]{Assumption}

\theoremstyle{definition}
\newtheorem{definition}[theorem]{Definition}
\newtheorem{remark}[theorem]{Remark}

\newcommand{\bgamma}{\mbox{\boldmath$\gamma$}}
\newcommand{\bvarphi}{\mbox{\boldmath$\varphi$}}
\newcommand{\balpha}{\mbox{\boldmath$\alpha$}}
\newcommand{\bbeta}{\mbox{\boldmath$\beta$}}
\newcommand{\bmu}{\mbox{\boldmath$\mu$}}
\newcommand{\sbeta}{{\mbox{\boldmath{\scriptsize{$\beta$}}}}}
\newcommand{\salpha}{{\mbox{\boldmath{\scriptsize{$\alpha$}}}}}
\newcommand{\sgamma}{{\mbox{\boldmath{\scriptsize{$\gamma$}}}}}
\newcommand{\szeta}{{\mbox{\boldmath{\scriptsize{$\szeta$}}}}}
\newcommand{\be}{\begin{equation}}
\newcommand{\ee}{\end{equation}}

\def\to{\longrightarrow}

\def \0{{\textbf{0}}}
\def\norm{{\| \kern -.05em | }}

\newcommand{\bc}{\begin{cases}}
\newcommand{\ec}{\end{cases}}

\title{Finite-dimensional representations\\ for controlled diffusions with delay}
\date{}
\author{Salvatore Federico\footnote{Universit\`a di Milano. E-mail: \texttt{salvatore.federico@unimi.it}} \and
  Peter Tankov\footnote{Universit\'e Paris Diderot. E-mail: \texttt{tankov@math.univ-paris-diderot.fr}}}  
\begin{document}
\maketitle
\abstract{
We study stochastic delay differential equations (SDDE) where the
coefficients depend on the moving averages of the state
process. As a first contribution, we provide sufficient
conditions under which a linear path functional of the solution of a SDDE admits a finite-dimensional Markovian
representation. As a second contribution, we show how approximate
finite-dimensional Markovian representations may be constructed when
these conditions are not satisfied, and provide an estimate of the
error corresponding to these approximations. These results are applied
to optimal control and optimal stopping problems for stochastic
systems with delay. 

\bigskip


\bigskip

}

\section{Introduction}{
In this paper we study a class of controlled stochastic differential equations
with memory, where the coefficients of the equation depend on the
moving average of the past values of the solution process  (so called distributed delay):
\begin{align}
dS_t=b\left(S_t,\int_{\mathbb{R}^-} \tilde{\alpha}_1(\xi)S_{t+\xi}d\xi, u_{t}\right)dt
+ \sigma\left(S_t,\int_{\mathbb{R}^-} \tilde{\beta}_1(\xi)S_{t+\xi}d\xi,u_{t}\right)dW_t,\label{intro.eq}
\end{align}
where $b,\sigma$,
$\tilde{\alpha}_1,\tilde{\beta}_1$ are given functions and $u=(u_{t})_{t\geq 0}$ is a control process. 
 Equations of this type appear in a variety of
domains such as economics \cite{GM06,GMS09} and finance
\cite{BBV10,BTW10,BV05,FP07,HR98}, as well as in physical sciences
\cite{M98}. In general this equation is  infinite-dimensional, which
means that it can be formulated as evolution equation in an
infinite-dimensional space of the form {$\mathbb{R}\times H_1$, where $H_1$ is a Hilbert space,} for the process $\mathbf X_t = (S_t, (S_{t+\xi})_{\xi\leq  0})$, but cannot be represented  via a
finite-dimensional controlled Markov process. This makes solving stochastic
control and optimal stopping problems associated to such systems
notoriously difficult. 

For this reason we are interested in
finding exact {- when possible -} or approximate finite dimensional representations for
functionals of solutions of \eqref{intro.eq}. We consider linear path
functionals of the form $Z_t = \langle \gamma, \mathbf X_t\rangle$,
where $\gamma$ is fixed. This functional may represent for example the
reward process of a control problem. {We say that the functional $Z$ admits a finite-dimensional Markovian representation if there exists a finite-dimensional subspace $V$ of the 
 space $\mathbb{R}\times H_1$ such that: 1) $V$ contains the vector $(1,0)\in \mathbb{R}\times H_1$; \ 2) the projection of the solution on this
subspace, call it $\mathbf X^V_t$, satisfies a finite dimensional stochastic
differential equation; 3) the functional $Z$ can be written as $Z_t = \langle \gamma, \mathbf
X^V_t\rangle$.} On the other hand, to find an approximate finite dimensional
representation for $Z$, we need to find a sequence of processes $(\mathbf
X^n_t)$ and a sequence of subspaces $(V_n)$, such that for every $n$,
the projection $\mathbf{X}^{n,V_n}$ satisfies a finite-dimensional SDE, and
such that for a sequence $(\gamma_n)$ to be determined,
$Z^n_t = \langle \gamma_n,
\mathbf X^{n,V}_t\rangle $ converges to $Z_t$ as $n\to \infty$. 

Our approach is different from most existing studies of invariance for
stochastic equations on Hilbert spaces (see e.g. \cite{Fil00,FT03}),
which require that the entire solution stays on a finite-dimensional
submanifold of the original space. Instead, we require that a
projection of the solution or an approximation thereof evolves on a finite-dimensional space. This projection only contains partial
information about the solution, but if the reward function of the
control problem only depends on this projection, this information is
sufficient to solve the control problem. 

Optimal control problems for stochastic systems with memory have been
considered by many authors starting with
\cite{kolmanovskii.maizenberg.73}. Solving the problem in the
infinite-dimensional setting being very difficult, some recent contributions
focus on special cases where the problem reduces to a
finite-dimensional one
\cite{bauer2005stochastic, EOS00, larssen2003hjb, oksendal2001maximum}. In
the general case, \cite{kushner2005numerical} extends the Markov chain
approximation method to stochastic equations with delay. A similar
method is developed in \cite{reiss.fischer.07}, and
\cite{fischer2008time} establish
convergence rates for an approximation of this kind. 
{
The infinite-dimensional Hilbertian approach to controlled
deterministic and stochastic systems with delays in the state variable
was employed in some papers.  For the deterministic case we can quote
\cite{FGG1, FGG2}, which perform a study of the
Hamilton-Jacobi-Bellman (HJB) equation in infinite
dimension\footnote{When the delay appears also in the control variable
  the infinite-dimensional representation is more involved. We refer
  to  \cite[Part\,II, Ch.\,4]{BDDM}, where a general theory is
  developed based on the paper \cite{VK}.}; for the stochastic case we
can quote \cite{F, GM06, GMS09} with some partial results on the
solution of the control problem (in \cite{GM06, GMS09} the delay is
considered also in the control variable, but the diffusion term is
just additive). {We should also mention the Banach space approach employed by \cite{FMT}:  the problem is embedded in the space of continuous functions and  the  HJB equation is approximated using the concept of mild solutions}.}
Optimal stopping problems for stochastic systems with
delay can be treated with methods similar to those used for optimal
control. \cite{federico2011optimal} and \cite{gapeev2006optimal}
discuss special cases where the infinite-dimensional problem reduces
to a finite-dimensional one. In the specific context of American
options written on the moving average of the asset price, \cite{BTW10}
propose a method based on Laguerre polynomial approximation, which is
extended and refined in the present paper. 

Let us now briefly summarize the contents of the paper. In section \ref{defs.sec} we
define the stochastic delay differential equation, state the
assumptions on the coefficients and introduce the main
notation. In Section \ref{inf.sec} we introduce and study an alternative
representation for this equation, as an evolution equation in an
infinite-dimensional Hilbert space. Section \ref{markov.sec} contains
the main results of the paper. First, we provide sufficient conditions for existence of an
exact finite dimensional Markovian representation for the output
process, namely that the coefficients belong to a certain
exponential-polynomial family (sum of exponential functions multiplied
by polynomials). Second, we describe a method for constructing an approximate finite-dimensional
representation using a specific exponential-polynomial family based on
Laguerre polynomials. The error of the approximation is also analyzed
here (Proposition \ref{err.prop}). Finally, Section \ref{app.sec} briefly discusses the
applications of our method to the solution of optimal control and
optimal stopping problems for stochastic systems with delay. Detailed
analysis of these applications and numerical examples is left for
further research. }
\section{The controlled stochastic delay differential equation }\label{defs.sec}
Let $(\Omega,\mathcal{F},\mathbb{P})$ be a complete probability space. On this space we suppose defined a Brownian motion $W=(W_t)_{t\geq 0}$ and denote by $\mathbb{F}=(\mathcal{F}_t)_{t\geq 0}$ the filtration generated by the process $W$ and enlarged by the $\mathbb{P}$-null sets.
\smallskip

Let $\mathbb{R}^-=(-\infty,0]$. To distinguish, deterministic functions  will be denoted with the time index in parentheses, while stochastic processes will be denoted with the time index as subscript.
\smallskip

Let $S=(S_t)_{t\geq 0}$ be a controlled diffusion on this space solving a stochastic delay differential equation (SDDE)
\be\label{SDDE1}
dS_t=b\left(S_t,\int_{\mathbb{R}^-} \tilde{\alpha}_1(\xi)S_{t+\xi}d\xi, u_{t}\right)dt\\
+ \sigma\left(S_t,\int_{\mathbb{R}^-} \tilde{\beta}_1(\xi)S_{t+\xi}d\xi,u_{t}\right)dW_t,
\ee
where $b,\sigma$,
$\tilde{\alpha}_1,\tilde{\beta}_1$ are given functions 
and $u=(u_{t})_{t\geq 0}$ is an adapted control process. 
Due to the dependence on the past, $S$ is not a controlled Markov diffusion. Moreover, in order to define the process $S$, 
one needs to
specify an initial condition not only at $t=0$, but also  for all $t<0$. In other terms, \eqref{SDDE}
has to be completed (in general) with an initial condition of the form
\be\label{IC}
S_0\ =\ s_0\ \in\ \mathbb{R}; \ \ \ \ \ \ \ S_{\xi}\ = \ s_1(\xi),\ \ \ \xi<0,
\ee
where $s_1$ is a given function. 
So, \emph{the initial datum is a function}. From what we have said
it is clear that, even if the process $S$ is one-dimensional, it may
not in general be represented as a finite-dimensional controlled Markov diffusion.\footnote{Nevertheless there are examples where a finite-dimensional Markovian representation can be obtained. We will study this kind of situation in Section \ref{finite.sec}, giving sufficient conditions for a finite-dimensional Markovian representation.
 } 
Moreover, even when the process $S$ is a controlled Markov diffusion,
i.e. $\alpha_1\equiv 0,\beta_1\equiv 0$ in \eqref{SDDE}, one may need to consider also the process\footnote{For example this process could appear in the cost functional of a control problem.}  
\be\label{function1}
 (Z_t)_{t\geq 0}\ =\ \left(\gamma_0S_t+\int_{\mathbb{R}^-}\tilde{\gamma}_1(\xi)S_{t+\xi}d\xi\right)_{t\geq 0},
\ee
where $\gamma_0\in\mathbb{R}$ and $\gamma_1$ is a function. Also in
this case it is clear that in general the Markovian representation of the system must be infinite-dimensional.
In this paper we deal  with the problem of  rewriting in an exact or approximate way the system above in terms of a Markov controlled finite-dimensional diffusion when at least one among the functions $\tilde{\alpha}_1,\tilde{\beta}_1,\tilde{\gamma}_1$ is not identically equal to $0$.
\smallskip
\begin{remark}
We stress that, altough we take  one-dimensional processes $W,S,Z$,
the argument can be easily generalized to the case of
multi-dimensional processes. Also we have taken an autonomous equation
for $S$, i.e.~there is no explicit time dependence of the coefficients
$b,\sigma$; this is done just for simplicity: all computations can be performed also in the non-autonomous case.
\end{remark}
{In the sequel, we are going to reformulate equation \eqref{SDDE1} as
an evolution equation in a Hilbert space. To allow a  set of initial data {possibly} containing the constant functions, we work with weighted spaces. }
We 
consider on $\mathbb{R}^-$ a weight function $w$ and make the
following standing assumption.
\begin{assumption}\label{ass:w}
$w\in C^{1}(\mathbb{R}^-;\mathbb{R})$, $w>0$, $w'/w$ is bounded.
\end{assumption}
 {Moreover, without any loss of generality, we also
   suppose that $w(0)=1$}.
Denote
\be\label{L2}
L^2_{w}\ \ :=\ \ L^2(\mathbb{R}^-,\, w(\xi) d\xi\,;\, \mathbb{R}).
\ee
When $w\equiv 1$ we simply denote the space above by $L^2$.
Throughout the paper, we shall work under the following assumptions on the model, guaranteeing existence and uniqueness for the solution to \eqref{SDDE1} and good properties for the problem we aim to study. 
\begin{assumption}\label{ass}
\begin{enumerate}
\item[]
\item The control process $u$ takes values in a set $U\subset \mathbb{R}^d$;
\item
$u\in\mathcal{U}$, where
$$\mathcal{U}\ \ =\ \ \{ (u_{t})_{t\geq 0} \ \mbox{adapted process belonging to }   L^2_{loc}(\mathbb{R}^+;L^2(\Omega;U))\};$$
\item
$b,\sigma:\mathbb{R}^2\times U\rightarrow \mathbb{R}$ are such that there exists constants $C_1,C_2\geq 0$ with
\begin{multline*}
|b(x,y,u)-b(x',y',u)|\ +\  |\sigma(x,y,u)-\sigma(x',y',u)|\\ \leq \ \ C_1\,(|x-x'|\ +\ |y-y'|), \ \ \  \ \forall x,x',y,y'\in\mathbb{R}, \ \forall u\in U;
\end{multline*}
$$|b(x,y,u)|\ +\  |\sigma(x,y,u)|\ \ \leq\ \  C_2\,(1+|x|+|y|), \ \ \ \ \ \ \ \ \  \forall x,y,\in\mathbb{R}, \ \forall u\in U;$$
\item {There exists $w$ satisfying Assumption \ref{ass:w} such that the functions $\tilde{\alpha}_1w^{-1/2},\tilde{\beta}_1w^{-1/2},\tilde{\gamma}_1w^{-1/2}$ belong to $L^2$.}
\end{enumerate}
\end{assumption}
\begin{remark}{
Typical weights are the exponential ones: $w(\xi)=e^{\lambda \xi}$,
$\lambda\in\mathbb{R}$. However, in some cases it may be necessary to
use other weight functions. For example, let $\tilde \gamma_1(\xi) =
\frac{1}{1+|\xi|^p}$ with $p>2$. Then, taking $w(\xi) =
\frac{1}{1+|\xi|^\frac{p}{2}}$, we ensure that simultaneously
$\tilde\gamma_1 w^{-1/2} \in L^2$ and the constant functions belong to
$L^2_w$. } {These two properties cannot hold simultaneously with an
exponential weight function.}
\end{remark}

\medskip
\noindent
{Define
$$
\alpha_1 \ =\ \tilde{\alpha}_1 w^{-1}\ \ \ \ \ \beta_1 \ =\ \tilde{\beta}_1 w^{-1} ,\ \ \ \ \ \gamma_1 \ =\ \tilde{\gamma}_1 w^{-1}.
$$}
Then, due to Assumption \ref{ass}(4),  
we have $\alpha_1,\beta_1,\gamma_1\in L^2_w$. Moreover,  \eqref{SDDE1}-\eqref{IC} can be rewritten as
\be
\begin{cases}
dS_t=b \left(S_t,\int_{\mathbb{R}^-} {\alpha}_1(\xi)S_{t+\xi}w(\xi)d\xi, u_{t}\right)dt + \sigma\!\left(S_t,\int_{\mathbb{R}^-} {\beta}_1(\xi)S_{t+\xi}w(\xi) d\xi,u_{t}\right)dW_t,\\\\
S_0=s_0, \ \ \ \ \  S_{\xi}=s_1(\xi), \ \ \xi<0,\label{SDDE}
\end{cases}
\ee
and \eqref{function1} can be rewritten as
\be\label{function}
Z\ =\ (Z_t)_{t\geq 0}\ =\ \left(\gamma_0S_t+\int_{\mathbb{R}^-}\gamma_1(\xi)S_{t+\xi}w(\xi)d\xi\right)_{t\geq 0},
\ee
\begin{proposition}
For every $\mathbf{s}:=({s}_0,s_1(\cdot))\in\mathbb{R}\times L^2_w$ and $u\in\mathcal{U}$,  \eqref{SDDE} admits a unique {up to indistinguishabiliy} strong solution $S^{\mathbf{s},u}$ on the probability space  $(\Omega,\mathcal{F},\mathbb{P})$, {and this solution admits a version with continuous paths.}
\end{proposition}
\textbf{Proof.} 
This result is an easy corollary of Theorem IX.2.1 in
\cite{revuzyor} {(note that without making any changes to the proof,
this theorem can be extended to the case when the coefficients depend
on a random adapted control).} 
{Let $T>0$ and define the maps $\tilde{B},\tilde{\Sigma} :[0,T]\times C([0,T];\mathbb{R})\times \Omega$
$$
\tilde{B}(t,{z}(\cdot),\omega) \ := \ b({z}(t), \int_{-\infty}^{-t} \alpha_1(\xi)s_1(t+\xi)
w(\xi) d\xi + \int_{-t}^0 z(t+\xi) \alpha_1(\xi) w(\xi) d\xi, u_t(\omega)),
$$
$$
\tilde{\Sigma} (t,z(\cdot),\omega) \ := \ b(z(t), \int_{-\infty}^{-t} \alpha_1(\xi)s_1(t+\xi)
w(\xi) d\xi + \int_{-t}^0 z(t+\xi)  \alpha_1(\xi) w(\xi) d\xi, u_t(\omega)).
$$
By  localizing in time, to  use the aforementioned result  we need to
check that 
\begin{enumerate}
\item
$ t\mapsto \tilde{B}(t,z(\cdot),\omega),  \ t\mapsto\tilde{\Sigma}(t,z(\cdot),\omega)$ are   bounded  for every given constant function $z(\cdot)\equiv z_0$ and uniformly in $\omega$;
 \item
  $\tilde{B}(t,\cdot,\omega), \tilde{\Sigma}(t,\cdot,\omega)$ are Lipschitz continuous, with respect to the uniform norm on $C([0,T];\mathbb{R})$,  uniformly in $t\in[0,T]$, $\omega\in \Omega$.
\end{enumerate}
}

{
Let us focus on $\tilde{B}$, the proof for $\tilde{\Sigma}$ being the same.
We first check the local boundedness on constant functions. Let $z(\cdot)\equiv z_0$. 
By the linear growth assumption,
\begin{align*}
|\tilde{B}(t,z(\cdot),\omega)| &= |b(z_0, \int_{-\infty}^{-t} \alpha_1(\xi)s_1(t+\xi)
w(\xi) d\xi + z_0\int_{-t}^0 \alpha_1(\xi) w(\xi) d\xi, u_t(\omega))| \\&\leq
C_2\left(1+ |y| + \left|\int_{-\infty}^{-t} \alpha_1(\xi)s_1(t+\xi)
w(\xi) d\xi\right| + |y|\left|\int_{-t}^0 \alpha_1(\xi) w(\xi) d\xi\right|\right).
\end{align*}
Let us denote the upper bound of $|w'|/w$ by $K$. Then, by Gronwall's
inequality, for all $a,b\leq 0$, 
\begin{align}
w(a)\ \leq  \ w(b) e^{K |b-a|}. \label{grw}
\end{align}
The term involving $s_1$ then satisfies
\begin{align*}
\left|\int_{-\infty}^{-t} \alpha_1(\xi)s_1(t+\xi)
w(\xi) d\xi\right|  &\leq \left(\int_{-\infty}^{-t} \alpha_1^2(\xi) w(\xi)
d\xi \right)^{\frac{1}{2}}\left(\int_{-\infty}^{-t} s_1^2(t+\xi) w(\xi)
d\xi \right)^{\frac{1}{2}}\\
&\leq \left(\int_{-\infty}^{0} \alpha_1^2(\xi) w(\xi)
d\xi \right)^{\frac{1}{2}}\left(\int_{-\infty}^{0} s_1^2(\xi) w(\xi-t)
d\xi \right)^{\frac{1}{2}}\\
&\leq e^{\frac{Kt}{2}} \left(\int_{-\infty}^{0} \alpha_1^2(\xi) w(\xi)
d\xi \right)^{\frac{1}{2}}\left(\int_{-\infty}^{0} s_1^2(\xi) w(\xi)
d\xi \right)^{\frac{1}{2}},
\end{align*}
which is  bounded in $[0,T]$ since $\alpha_1, s_1\in L^2_w$. Similarly,
\begin{align}
\left|\int_{-t}^0 \alpha_1(\xi) w(\xi) d\xi\right|\ \leq \ \int_{-t}^0 |\alpha_1(\xi)| w(\xi) d\xi\ \leq\ \left(\int_{-t}^0
\alpha_1^2(\xi) w(\xi) d\xi \right)^{\frac{1}{2}}\left(\int_{-t}^0
 w(\xi) d\xi \right)^{\frac{1}{2}},\label{csalpha}
\end{align}
which is  bounded in $[0,T]$ as well since $\alpha_1\in L^2_w$ and $w\in C^1(\mathbb{R}^-;\mathbb{R})$.}

{
Let us now check the  Lipschitz condition.  By the Lipschitz
property of $b$, 
\begin{align*}
|B(t,z(\cdot), \omega) - B(t,z'(\cdot),\omega)| &\leq C_1\Big(|z(t) - z'(t)| + \int_{-t}^0
|z({t+\xi}) - z'({t+\xi})| |\alpha_1(\xi)| w(\xi) d\xi\Big)\\
&\leq C_1\Big(1 +\int_{-T}^0 |\alpha_1(\xi)| w(\xi) d\xi \Big)\cdot \max_{0\leq u \leq T} |z(u) - z'(u)|,
\end{align*}
and the functional Lipschitz property follows from inequality \eqref{csalpha}. }
\hfill$\square$

\section{Product space infinite-dimensional representation}\label{inf.sec}
In this section we provide an infinite-dimensional representation of SDDE \eqref{SDDE} in the product Hilbert space
$$H_w\ \ :=\ \ \mathbb{R}\times L^2_w.$$
When $w\equiv 1$ we simply denote the space above by $H$.
We denote by $\mathbf{x}=(x_0,x_1)$ the generic element of $H_w$,
noting that the second component is a function. The norm and the inner
product of $H_w$, defined in the usual way from the norm and the inner products of the two components, will be denoted, respectively,  by $\|\cdot\|_w$,  $\langle\cdot,\cdot\rangle_w$. Also, when $w\equiv 1$ we simply denote the norm and the inner product above by $\|\cdot\|$,  $\langle\cdot,\cdot\rangle$.
\subsection{Preliminaries}
{
Let us introduce the weighted Sobolev spaces on $\mathbb{R}^-$ as follows (we refer to \cite[Ch.\,VIII]{B} for an introduction to Sobolev spaces under the Lebesgue measure on intervals).  Given $f\in L^2_{loc}(\mathbb{R}^-;\mathbb{R})$, we say that $f$ admits weak derivative $g\in L^2_{loc}(\mathbb{R}^-;\mathbb{R})$ if  
$$ \int_{\mathbb{R}} f(\xi)\varphi'(\xi)d\xi \ \ = \ \ - \int_{\mathbb{R}} g(\xi)\varphi(\xi)d\xi,  \ \ \ \ \forall \varphi\in C_c^1((-\infty,0);\mathbb{R}).$$
It is well known that, if such a function $g$ exists, it is
unique. Moreover it coincides with the classical derivative $f'$ when
$f\in C^1(\mathbb{R}^-;\mathbb{R})$. By extension, the function $g$ is denoted by $f'$ in general, i.e.,~also when $f\notin C^1(\mathbb{R}^-;\mathbb{R})$.  We denote the space of functions of $L^2_{loc}(\mathbb{R}^-;\mathbb{R})$  admitting weak derivative in $L^2_{loc}(\mathbb{R}^-;\mathbb{R})$  by $W^{1,2}_{loc}(\mathbb{R}^-;\mathbb{R})$.
%
 \cite[Th.\,VIII.2]{B} states that for every $f\in W^{1,2}_{loc}(\mathbb{R}^-;\mathbb{R})$ there
 exists a locally absolutely continuous version of $f$ on $\mathbb{R}^-$, so that it holds
\be\label{abscont}
f(\xi)-f(\xi_0)\ \ =\ \ \int_{\xi_0}^{\xi} f'(r)dr, \ \ \ \  \forall \ \xi_0 \leq  \xi\leq 0.
\ee
Given $f\in W^{1,2}_{loc}(\mathbb{R}^-;\mathbb{R})$, we shall always refer to its absolutely continuous version. By
\cite[Cor.\,VIII.10]{B}, if $f,g\in W^{1,2}_{loc}(\mathbb{R}^-;\mathbb{R})$, then   $fg\in W^{1,2}_{loc}(\mathbb{R}^-;\mathbb{R})$ and
\be\label{dpr}
(fg)'\ \ =\ \ f'g+fg',
\ee
so the integration by parts formula
$$
\int_a^b f'(\xi)g(\xi) d\xi \ \ =\ \ f(b)g(b)-f(a)g(a) - \int_{a}^b f(\xi)g'(\xi)d\xi, \ \ \ \ \forall \ a\leq b\leq 0,
$$
holds true for all $f,g\in W^{1,2}_{loc}(\mathbb{R}^-;\mathbb{R})$.
On the elements of the space $W^{1,2}_{loc}(\mathbb{R}^-;\mathbb{R})$ we define the norm
$$
\|f\|_{W^{1,2}_w} \ :=\ \int_{\mathbb{R}^-} (|f(\xi)|^2+|f'(\xi)|^2) w(\xi)d\xi,
$$ 
and, moreover, we define  the space 
$$
W^{1,2}_w \ \ := \ \ \left\{ f \in W^{1,2}_{loc} \ |  \ \|f\|_{W^{1,2}_w}<\infty\right\}.
$$
Clearly $W^{1,2}_w\subset L^2_w$.
 The linear maps
$$
\begin{array}{cccc}
 & (W^{1,2}_w,\|\cdot\|_{W^{1,2}_w}) & \longrightarrow & (L^2_w,\|\cdot\|_{L^2_w}) \times (L^2_w,\|\cdot\|_{L^2_w}) ,\\
& f&\longmapsto & (f,f'),
\end{array}
$$
and 
$$
\begin{array}{cccc}
 & (L^2_w,\|\cdot\|_{L^2_w}) & \longrightarrow & (L^2,\|\cdot\|_{L^2}) ,\\
& f&\longmapsto & fw^{\frac{1}{2}},
\end{array}
$$
are isometries, so, since $L^2$ is a separable Banach space,   we deduce that  $W^{1,2}_w$ is a separable Hilbert space when endowed with the inner product
$$
\langle f,g\rangle_{W^{1,2}_w} \ :=\ \int_{\mathbb{R}^-} f(\xi)g(\xi)w(\xi)d\xi. 
$$ 
{By the assumption $w'/w$ bounded, 
denoting 
 the upper bound of $|w'|/w$ by $K$, we see that if $f\in W^{1,2}_w,$ then 
 \begin{align}\label{estt}
\|f{w}^{1/2}\|^2_{W^{1,2}} &= \int_{\mathbb R^-} \Bigg\{f^2(\xi) w(\xi) +
\Big(f'(\xi)\sqrt{w(\xi)} +
f(\xi)\frac{w'(\xi)}{2\sqrt{w(\xi)}}\Big)^2\Bigg\}d\xi\nonumber\\ &\leq
\int_{\mathbb R^-}(f^2(\xi) + 2f'(\xi)^2 + \frac{K^2}{2}
f^2(\xi))w(\xi) d\xi \leq (2+\frac{K^2}{2}) \|f\|_{W^{1,2}_w}^2.
\end{align}
Thus, if $f\in W^{1,2}_w,$ then 
  then $fw^{1/2}\in W^{1,2}$.} Hence,  Corollary \cite[Cor.\,VIII.8]{B} applied to $fw^{1/2}$ yields
$$
\lim_{\xi\rightarrow -\infty} f(\xi)w^{\frac{1}{2}}(\xi) \ =\ 0, \ \ \ \forall f\in W^{1,2}_w.
$$
Recalling again our assumptions on $w$,  we see that the following weighted  integration by parts formula holds for all $f,g\in W^{1,2}_w$:
\begin{equation}\label{IBP}
\int_{\mathbb{R}^-} f'(\xi)g(\xi) w(\xi)d\xi \ \ =\ \ f(0)g(0)-\int_{\mathbb{R}^-} f(\xi)\bigg(g'(\xi)+g(\xi)\frac{w'(\xi)}{w(\xi)}\bigg)w(\xi) d\xi.
\end{equation}
Now, consider on the space $H_w$ the family of linear bounded operators $(T(t))_{t\geq 0}$   acting as follows: 
\be\label{expl:sem}
T(t)\mathbf{x}\ =\ \left([T(t)\mathbf{x}]_0,\,[T(t)\mathbf{x}]_1\right) \ =\ \left(x_0, \ x_0\mathbf{1}_{{(0,t]}}(t+\cdot)+x_1(t+\cdot)\mathbf{1}_{\mathbb{R}^-}(t+\cdot)\right).
\ee
Simple computations show that
\be\label{estsem}
\|T(t)\|_{\mathcal{L}(H_w)}\ \ \leq \ \ 1+t, \ \ \ \ \ \forall  t\geq 0.
\ee
We are going to study the semigroup properties of $(T(t))_{t\geq 0}$. For basic facts about the theory of semigroups we refer to the classical monographs \cite{Davies, EN}.
\begin{proposition}\label{prop:semi}
The family of linear operators $(T(t))_{t\geq 0}$ defined in \eqref{expl:sem} is a strongly continuous  semigroup on the space $H_w$, generated by the closed unbounded operator  
$\mathcal{A}$ defined on
\be\label{DA}
D(\mathcal{A})\ \ =\ \ \{\mathbf{x}=(x_0,x_1)\in H_w \ | \ x_1\in W^{1,2}_w, \ x_0=x_1(0) \}
\ee
by  
\be\label{AA}
\mathcal{A}\mathbf{x}\ =\   \left(0, \ x_1'\right).
\ee
\end{proposition}
\textbf{Proof.} 
The fact that $(T(t))_{t\geq 0}$ is a semigroup is immediate by the definition. The fact that it is strongly continuous follows by the continuity of translations in $L^2_w$, { which can be proved, e.g., starting from the continuity of translation in $L^2$ and exploiting \eqref{grw}.}

Now let us show that $(T(t))_{t\geq 0}$ is generated by $\mathcal{A}$.  Set
$$
\mathcal{D} \ \ := \ \ \{\mathbf{x}=(x_0,x_1)\in H_w \ | \ x_1\in W^{1,2}_w, \ x_0=x_1(0) \}
$$
and 
take $\mathbf{x}\in\mathcal{D}$. Since $x_1\in W^{1,2}_w$, it is
absolutely continuous. So, extending $x_1$ to $\mathbb{R}$ by setting $x_1(\xi)=x_1(0)$ for $\xi>0$,
we can write 
\begin{eqnarray*}
x_1(t+\xi)-x_1(\xi)&=& t \int_0^1x_1'(\xi+\lambda t) d\lambda, \ \ \ \ \forall \xi\in \mathbb{R}^-, \ \forall t\geq 0.
\end{eqnarray*}
Hence, taking into account that $x_1(0)=x_0$, we have
\begin{eqnarray*}
\left\|\frac{T(t)\mathbf{x}-\mathbf{x}}{t}-(0,x_1')\right\|^2_{w} & = & \int_{\mathbb{R}^-}w(\xi)d\xi\left|\int_0^1(x_1'(\xi+\lambda t)-x_1'(\xi)) d\lambda\right|^2\\
&\leq & \int_0^1d\lambda \int_{\mathbb{R}^-} \left|x_1'(\xi+\lambda t)-x_1'(\xi)\right|^2w(\xi)d\xi\\
&=&  \int_0^1d\lambda \int_{\mathbb{R}^-} \left\| T(\lambda t)(0,x_1')-(0,x_1')\right\|^2_{H_w}d\xi.
\end{eqnarray*}
By \eqref{estsem} and  from the inequality above, we can conclude by dominated convergence that $\mathcal{D}\subset D(\mathcal{A})$ and that $\mathcal{A}$ acts as stated in \eqref{AA} on the elements of $\mathcal{D}$. 

We need now to show that $\mathcal{D}=D(\mathcal{A})$. For that, we
notice that $\mathcal{D}$ is clearly dense in $H_w$ and that
$T(t)\mathcal{D}\subset \mathcal{D}$  for any $t\geq 0$. Hence, {  by
\cite[Ch.\,II, Prop.\,1.7, p.\,53]{EN}, $\mathcal{D}$ is a core for $D(\mathcal{A})$ (i.e. is dense in $D(\mathcal{A})$ endowed with the graph norm  $\|\cdot\|_{D(\mathcal{A})}$).} Hence, it just remains to show that $\mathcal{D}$ is closed with respect to the graph norm to conclude $\mathcal{D}= D(\mathcal{A})$. So,   take a sequence $(\mathbf{x}^n)=(x^n_0,x^n_1)\subset \mathcal{D}$ converging with respect to $\|\cdot\|_{D(\mathcal{A})}$ to some $\mathbf{x}=(x_0,x_1)\in H_w$. Then we have
\begin{equation}\label{convvv}
x^n_0\rightarrow x_0\ \ \  \mbox{in}\ \mathbb{R}; \ \ \ \ \ \ \ x^n_1\rightarrow x_1 \ \ \  \mbox{in}\  W^{1,2}_w.
\end{equation}
{
We immediately  deduce that $x_1\in W^{1,2}_w$. 
By \eqref{estt}, the linear map
$$
\mathcal{L}:\  W^{1,2}_w\ \rightarrow \ W^{1,2}, \ \ \ f\ \mapsto\ fw^{1/2}, 
$$
is continuous.
Since we have (see, e.g. \cite[Th.\,8.8]{B}) the Sobolev continuous embedding $\iota: \ W^{1,2}\hookrightarrow L^\infty(\mathbb{R}^-;\mathbb{R})$, the map $\iota\circ \mathcal{L}$ is continuous. Taking into account also that $x_1$ is absolutely continuous,  
we  deduce from the second convergence in \eqref{convvv}
$$x^n_1(0)\rightarrow x_1(0) \ \ \ \mbox{in} \ \mathbb{R}.$$
Since $x^n_0=x^n_1(0)$, we conclude $x_1(0)=x_0$, and the proof is complete.}
\hfill$\square$
\subsection{Infinite-dimensional representation} 
Define the elements of $H_{w}$
$$\mathbf{e}^0\ :=\ (1,0),\ \ \ \ \ \balpha\ :=\ (0,\alpha_1),\ \ \  \ \ \bbeta\ :=\  (0,\beta_1),$$
and
the Lipschitz continuous nonlinear operators
$$B:\ H_w\times U\ \rightarrow \ H_w, \ \ \ \ \  \ B(\mathbf{x},u)\ :=\  \left(b(\langle\mathbf{e}^0, \mathbf{x}\rangle_w, \langle\balpha, \mathbf{x}\rangle_w,u),0\right);$$  
$$ \Sigma:\ H_w \times U\ \rightarrow \ H_w, \ \ \ \ \ \Sigma(\mathbf{x},u)\ :=\  \left(\sigma(\langle\mathbf{e}^0, \mathbf{x}\rangle_w, \langle\bbeta, \mathbf{x}\rangle_w,u),0\right).$$ 
Given $\mathbf{x}\in H_w$ and $u\in\mathcal{U}$, we consider the following  stochastic evolution equation in the space $H_w$:
\be\label{eq:inf}
\bc
d\mathbf{X}_{t}\ \ =\ \ \mathcal{A}\mathbf{X}_{t}dt+B(\mathbf{X}_{t},u_{t})dt+ \Sigma(\mathbf{X}_{t},u_{t}) dW_t,\\
\mathbf{X}_{0}\ \ =\ \ \mathbf{x},
\ec
\ee
At least formally \eqref{eq:inf} should represent \eqref{SDDE} in  $H_w$: if there exists a  unique solution (in some sense) $\mathbf{X}$  to \eqref{eq:inf}, we expect that 
$$\mathbf{X}_{t}\ \ =\ \ \left(S_t,(S_{t+\xi})_{\xi\in\mathbb{R}^-}\right), \ \ \ \ \ \ \ \forall t\geq 0,$$
where $S$ is the solution to \eqref{SDDE} with $\mathbf{s}=\mathbf{x}$.
{We notice that \eqref{eq:inf} is an equation in infinite dimension, but the noise is one-dimensional\footnote{{In the usual language of stochastic integration in infinite-dimension (see \cite{DPZ, GM11, PR}), $\Sigma(t,\mathbf{x})$ can be seen as a Hilbert-Schmidt operator from $\mathbb{R}$ to $H$. 
}}.}

We are going to introduce two concepts of solution to \eqref{eq:inf},
which in this case coincide with each other. Before we introduce the  operator $\mathcal{A}^*$ adjoint of $\mathcal{A}$. 
\begin{proposition} The adjoint $\mathcal{A}^*$ of the operator $\mathcal{A}$ is defined on 
\be\label{D(A*)}
{D}(\mathcal{A}^*)\ \ =\ \ \{ \mathbf{x}\in H_{w} \ | \ x_1\in W^{1,2}_{w}\},
\ee
by
\be\label{{A}^*}
\mathcal{A}^*\mathbf{x}\ \ =\ \  \left(x_1(0),\,-x_1'-x_1\frac{w'}{w}\right).
\ee
\end{proposition}
\textbf{Proof.} Let $\mathcal{D}:=  \{ \mathbf{x}\in H_{w} \ | \ x_1\in W^{1,2}_{w}\}$, $\mathbf{x}\in \mathcal{D}$ and $\mathbf{y}\in D(\mathcal{A})$. Using  \eqref{IBP} and the fact that $y_1(0)=y_0$, we can write
\begin{eqnarray*}
\langle \mathcal{A}\mathbf{y}, \mathbf{x}\rangle_w&=& \int_{\mathbb{R}^-} y'_1(\xi)x_1(\xi)w(\xi) d\xi\\
&=&y_1(0)x_1(0)-\int_{\mathbb{R}^-} y_1(\xi)\left(x_1'(\xi)+x_1(\xi)\frac{w'(\xi)}{w(\xi)}\right)w(\xi)d\xi\\
&=&y_0x_1(0)-\int_{\mathbb{R}^-} y_1(\xi)\left(x_1'(\xi)+x_1(\xi)\frac{w'(\xi)}{w(\xi)}\right)w(\xi)d\xi.
\end{eqnarray*}
So, we can conclude  that $\mathcal{D}\subset D(\mathcal{A}^*)$, and that $\mathcal{A}^*$ acts as in \eqref{{A}^*} on $\mathcal{D}$.

{
Now let us show that actually $D(\mathcal{A}^*)=\mathcal{D}$. Simple computations shows that the expression of the adjoint semigroup of $T(\cdot)$ in $H_w$ is 
$$
T^*(t)\ =\ \left(x_0+\int_{-t}^0x_1(\xi), \ x_1(\cdot-t)\frac{w(\cdot-t)}{w(\cdot)}\right),\ \ \ \ \mathbf{x}\in H_w.
$$
The set $\mathcal{D}$ is clearly dense in $H_w$  and 
$T^*(t)\mathcal{D}\subset \mathcal{D}$  for any $t\geq 0$. Hence,   by
\cite[Ch.\,II, Prop.\,1.7, p.\,53]{EN}, $\mathcal{D}$ is a core for $D(\mathcal{A}^*)$. On the other hand, in analogy with the proof of Proposition \ref{prop:semi}, one can show that $\mathcal{D}$ is closed with respect to the graph norm, so we  conclude that $\mathcal{D}= D(\mathcal{A}^*)$.}
\hfill$\square$
\begin{definition}\label{def:mild-weak}
\begin{itemize}
\item[(i)]
Let $\mathbf{x}\in H_w$, and let  $u\in\mathcal{U}$. An adapted process $\mathbf{X}=\mathbf{X}^{\mathbf{x},u}\in L^2_{loc}(\mathbb{R}^+;L^2(\Omega;H_w))$ is called  mild solution to \eqref{eq:inf} if for every $t\geq 0$
\begin{align}
\mathbf{X}_{t}\ \ = \ \ T(t)\mathbf{x}\ +\ \int_0^tT(t-r)B(\mathbf{X}_{r},u(r))  dr \ +\  \int_0^t T(t-r)  \Sigma (\mathbf{X}_{r},u(r))dW_r.\label{eq:mild}
\end{align}
\item[(ii)]
Let $\mathbf{x}\in H_w$, $u\in\mathcal{U}$.
An adapted process $\mathbf{X}=\mathbf{X}^{\mathbf{x},u}\in L^2_{loc}(\mathbb{R}^+;L^2(\Omega;H_w))$  is called  weak solution to \eqref{eq:inf} if for each $\bvarphi\in{D}(\mathcal{A}^*)$ and every $t\geq 0$
\begin{multline*}
\langle \mathbf{X}_{t},\bvarphi\rangle_w\ =\ \langle \mathbf{x},\bvarphi\rangle_w \,+\,\int_0^t\langle \mathbf{X}_{s}, \mathcal{A}^*\bvarphi\rangle_w ds \\+\, \int_0^t\langle B(\mathbf{X}_{s},u_{s}),\bvarphi\rangle_w ds\,+\,
\int_0^t \langle \Sigma(\mathbf{X}_{s},u_{s}),  \bvarphi\rangle_w dW_s.
\end{multline*}
\end{itemize}
\end{definition}
\begin{theorem}
For each $\mathbf{x}\in H_w$ and $u\in\mathcal{U}$,  the SDE \eqref{eq:inf} admits a unique  { up to indistinguishabiliy} mild solution $\mathbf{X}=\mathbf{X}^{\mathbf{x},u}$ which coincides with the unique weak solution, {and this solution admits a version with continuous paths}. 

Moreover, we have the equality in $L^2(\Omega,\mathcal{F},\mathbb{P};H_w)$
$$\mathbf{X}_{t}\ \ =\ \ \left(S_t,\,(S_{t+\xi})_{\xi\in\mathbb{R}^-}\right), \ \ \ \ \ \forall t\geq 0,$$
where $S$ is the solution to \eqref{SDDE} under the control $u$ and with initial datum $\mathbf{s}=\mathbf{x}$.
\end{theorem}
\textbf{Proof.} The existence and uniqueness of a mild solution, as well as the fact that it coincides with the (unique) weal solution,   is a straightforward application of the theory of infinite-dimensional stochastic differential equations (see, e.g.,  \cite[Ch.\,7]{DPZ}) and \cite[Ch.\,3]{GM11}). 

{ For the second part of the claim, let $S$ be the solution of \eqref{SDDE} and
define
$$\widetilde{\mathbf{X}}_{t}\ :=\ \left(S_t,\,(S_{t+\xi})_{\xi\in\mathbb{R}^-}\right), \ \ \ \mathbf x\ :=\ (s_0, (s_1(\xi))_{\xi\in \mathbb R^-}).$$ Then, 
\begin{align*}
S_t &= s_0 + \int_0^t b(S_r, \int_{\mathbb R^-} S_{r+\xi} \alpha_1(\xi)
w(\xi )d\xi , u(r))dr + \int_0^t \sigma(S_r, \int_{\mathbb R^-}
S_{r+\xi} \beta_1(\xi) w(\xi )d\xi, u(r))dW_r\\
& = s_0 + \int_0^t b(\langle \mathbf e^0,\widetilde{\mathbf X}_r\rangle_w, \langle \mathbf
\alpha, \widetilde{\mathbf X}_r\rangle_w, u(r))dr + \int_0^t \sigma(\langle \mathbf e^0,\widetilde{\mathbf X}_r\rangle_w, \langle \mathbf
\beta, \widetilde{\mathbf X}_r\rangle_w, u(r))dW_r \\
&=[T(t)\mathbf x]_0 + \int_0^t [T(t-r) B(\widetilde{\mathbf X}_r,u(r))]_0 dr
+  \int_0^t [T(t-r) \Sigma(\widetilde{\mathbf X}_r,u(r))]_0 dW_r
\end{align*}
and for every $\xi \in \mathbb R^-$,
\begin{align*}
S_{t+\xi} & = \mathbf 1_{t+\xi
  <0}  s_1(t+\xi) + \mathbf 1_{t+\xi
  \geq 0}\Bigg\{ s_0 + \int_0^{t+\xi}  b(S_r, \int_{\mathbb R^-} S_{r+\eta} \alpha_1(\eta)
w(\eta )d\eta , u(r))dr\\ &\qquad \qquad \qquad \qquad +\int_0^{t+\xi} \sigma(S_r, \int_{\mathbb R^-}
S_{r+\eta} \beta_1(\eta) w(\eta )d\eta, u(r))dW_r
\Bigg\}\\
& = \mathbf 1_{t+\xi
  <0}  s_1(t+\xi) + \mathbf 1_{t+\xi
  \geq 0}\Bigg\{ s_0 + \int_0^{t+\xi}  b(\langle \mathbf e^0,\widetilde{\mathbf X}_r\rangle_w, \langle \mathbf
\alpha, \widetilde{\mathbf X}_r\rangle_w, u(r))dr \\ & \qquad \qquad \qquad \qquad+\int_0^{t+\xi} \sigma(\langle \mathbf e^0,\widetilde{\mathbf X}_r\rangle_w, \langle \mathbf
\beta, \widetilde{\mathbf X}_r\rangle_w, u(r))dW_r
\Bigg\}\\
& = s_0 \mathbf 1_{t+\xi \geq 0} + s_1(t+\xi)\mathbf 1_{t+\xi
  <0}+\int_0^t \mathbf 1_{t-r+\xi \geq 0} \,b(\langle \mathbf e^0,\widetilde{\mathbf X}_r\rangle_w, \langle \mathbf
\alpha, \widetilde{\mathbf X}_r\rangle_w, u(r))dr\\&\qquad \qquad \qquad \qquad+\int_0^t \mathbf 1_{t-r+\xi \geq 0} \,\sigma(\langle \mathbf e^0,\widetilde{\mathbf X}_r\rangle_w, \langle \mathbf
\beta, \widetilde{\mathbf X}_r\rangle_w, u(r))dW_r\\
&=\left\{[T(t)\mathbf x]_1 + \int_0^t [T(t-r) B(\widetilde{\mathbf X}_r,u(r))]_1 dr
+  \int_0^t [T(t-r) \Sigma(\widetilde{\mathbf X}_r,u(r))]_1 dW_r\right\}\Bigg|_{\xi}
\end{align*}
which shows that $\widetilde{\mathbf X}$ satisfies \eqref{eq:mild} and
therefore coincides with the unique mild solution. }
\hfill$\square$\\\\
Since the two concepts of solutions coincide each other in this case, from now on we just say solution to refer to the mild or weak solution.
The following technical result will be used in the following section. 
\begin{proposition}\label{prop:growthX}
Let $\mathbf{X}=\mathbf{X}^{\mathbf{x},u}$ be the solution to \eqref{eq:inf}. Then
{$$\mathbb{E}\left[\sup_{0\leq t \leq T}\|\mathbf{X}_{t}\|_w^2\right]\ \ \leq\ \  p_1(T)\|\mathbf{x}\|_w^{2}+p_2(T),\ \ \ \ \ \ \ \forall T\geq 0$$}
where $p_1,p_2$ are {locally bounded functions.}
\end{proposition}
\textbf{Proof.} 
{
We notice that
\be\label{estB}
\|B(\mathbf{x},u)\|_w\ \ \leq \ \ C_{b,\salpha}(1+\|\mathbf{x}\|_w), \ \ \ \ \  \|\Sigma(\mathbf{x},u)\|_w\ \ \leq \ \ C_{\sigma,\sbeta}(1+\|\mathbf{x}\|_w),
\ee
 where 
 $$
 C_{b,\salpha}\ \ =\ \ C_2(1+\|\balpha\|_w),\ \ \ \  C_{\sigma, \sbeta}\ \ =\ \ C_2(1+\|\bbeta\|_w).
$$}
Let $T>0$. Using Definition \ref{def:mild-weak}-(i) and \eqref{estsem}, we have
\begin{align}
\mathbb{E}\left[\sup_{0\leq t \leq T}\|\mathbf{X}_{t}\|_w^2\right] &\leq
3(1+T)^2\|\mathbf{x}\|_w^2 + 3 \mathbb{E}\left[\sup_{0\leq t \leq T}\Bigg\|\int_0^t T(t-r) B(\mathbf
X_r,u(r))dr\Bigg\|_w^2\right] \notag\\&+ 3 \mathbb{E}\left[\sup_{0\leq t \leq T}\Bigg\|\int_0^t T(t-r) \Sigma(\mathbf X_r,u(r))dW_r\Bigg\|_w^2\right].\label{3term}
\end{align}
By Jensen's inequality (using the convexity of $\|\cdot\|_w^2$) and by the estimates \eqref{estB} and \eqref{estsem},
we deduce
\begin{align}
\Bigg\|\int_0^t T(t-r) B(\mathbf
X_r,u(r))dr\Bigg\|_w^2 &\leq t \int_0^t \|T(t-r)B(\mathbf
X_r,u(r)) \|_w^2 dr\notag
\\
&\leq t\int_0^t (1+t-r)^2
C^2_{b,\salpha}(1+\|\mathbf{X}_r\|_w)^2dr\notag\\
& \leq {2T(1+T)^2 C^2_{b,\salpha} \Big(\int_0^T(1+
\|\mathbf{X}_r\|_w^2) dr\Big).} \label{convb}
\end{align}

{
On the other hand,  the estimates \eqref{estB} and \eqref{estsem},  Doob's inequality and It\^o's isometry in infinite dimension (see, e.g., \cite[Ch,\,2]{GM11} or \cite[Ch.\,2]{PR}) also yield
\begin{align*}
&\mathbb E\left[\sup_{0\leq t \leq T}\Bigg\|\int_0^t T(t-r) \Sigma(\mathbf
X_r,u(r))dW_r\Bigg\|_w^2\right] \ \leq \ 8 (1+T)^2C_{\sigma,\bbeta}^2\int_0^T
(1+\mathbb{E}\|\mathbf X_r\|^2_w)dr.
\end{align*}
The claim follows from Gronwall's inequality.
}
\hfill$\square$
\section{Markovian  representations}\label{markov.sec}
{In this section we give conditions for the existence of exact representations
and provide a method to construct approximate representations for the process $Z$ as a deterministic function of the current state of a 
finite-dimensional controlled Markov diffusion.}
\subsection{Preliminaries}
{The first step is to characterize the finite-dimensional
  subspaces of $H_w$ which are stable with respect to the operator
  $\mathcal A^*$, which generates the infinite-dimensional
  structure of our delay equation.
 }

Set
$$\lambda^*\ \ := \ \ \inf\,\{\lambda\in \mathbb{R} \ | \ e^{\lambda\xi}w^{-\frac{1}{2}}(\xi)\in L^2\}.$$
We assume, just for simplicity, that the infimum above is not attained, i.e. that $\int_{\mathbb{R}^-} e^{2\lambda^*\xi} w^{-1}(\xi)d\xi=\infty.$ 
Let us introduce the differential operator
$$
\mathcal{D}_w: \ W_w^{1,2} \ \longrightarrow \ H_w, \ \ \ \ \ \ v \ \longmapsto\  -\frac{(vw)'}{w}.
$$
\begin{definition}
We say that a finite dimensional subspace $\mathcal{V}$ of $L^2_w$ is $\mathcal{D}_w$-stable if $\mathcal{V}\subset W^{1,2}_w$ and $\mathcal{D}_w \mathcal{V}\subset \mathcal{V}$.
\end{definition}
We have the following characterization of $\mathcal{D}_w$-stable subspace.
\begin{proposition}\label{prop:stable}
$\mathcal{V}$ is an $n$-dimensional $\mathcal{D}_w$-stable subspace of $L^2_w$ if and only if 
\be\label{Vstable}
\mathcal{V}\ \ =\ \  \mbox{\mbox{Span}}\ \{v_1w^{-1}, ..., v_n w^{-1}\},
\ee 
where  $\{v_1,...,v_n\}$ is a set of linearly independent functions such that the vector function $\mathbf{v}= (v_1,..., v_n)$ is 
  solution of the vector-valued ODE
 \be\label{Mstable}
 \mathbf{v}'(\xi)\ \ =\ \ \mathbf{M} \mathbf{v}(\xi), \ \ \ \xi\leq 0, 
 \ee 
 for some matrix $\mathbf{M}\in\mathbb{R}^{n\times n}$ whose
 eigenvalues have real part strictly greater than $\lambda^*$.
 \end{proposition}
\textbf{Proof.} Let $\mathcal{V}$ be in the form \eqref{Vstable} with  $\{v_1,...,v_n\}$ linearly independent functions such that $\textbf{v}=(v_1,...,v_n)$ is  solution of \eqref{Mstable} for some matrix $\textbf{M}\in\mathbb{R}^{n\times n}$ with eigenvalues $(\lambda_p)_{p=1,...,k}$ such that $Re(\lambda_p)>\lambda^*$ for all $p=1,...,k$.
Clearly $\mbox{dim}(V)=n$. Since $v$ solves  \eqref{Mstable}, the $v^i$'s are linear combination of functions of the form 
\be\label{str}
\xi^j e^{Re(\lambda_p)\xi}\cos(Im(\lambda_p)\xi), \ \ \ \ \xi^j e^{Re(\lambda_p)\xi}\sin(Im(\lambda_p)\xi).
\ee
Then,
 since  $Re(\lambda_p)>\lambda^*$ for all $p=1,...,k$, by definition of $\lambda^*$ we see that 
 $$
w^{-1}v_i\ \ \in\ \ W^{1,2}_w, \ \ \ \ \forall i=1,...,n,
 $$
 and therefore
 $\mathcal{V}\,\subset\, W^{1,2}_w$. On the other hand, given $f\in \mathcal{V}$, we have 
 $$f\ \  =\ \ w^{-1}\sum_{i=1}^n \mu_i v_i\ \ =\ \ w^{-1} \bmu^T \mathbf{v}, \ \ \ \mbox{ for some } \ {\bmu}=(\mu_i)_{i=1,...,n}\in\mathbb{R}^n.$$ 
Hence, since  $\mathbf{v}$ solves \eqref{Mstable}, we see that   
$$
\mathcal{D}_w f\ \ =\ \ -\frac{(w w^{-1}{\bmu}^T\textbf{v})'}{w}\ \ =\ \ -w^{-1} {\bmu}^T \mathbf{v}'\ \ =\ \ -w^{-1}{\bmu}^T\mathbf{M}\mathbf{v}\ \ \in\ \ \mathcal{V},
$$
showing the $\mathcal{D}_w$-stability of $\mathcal{V}$.

Conversely, let us suppose that  $\mathcal{V}$ is an $n$-dimensional $\mathcal{D}_w$-stable subspace, and let $\{\tilde{v}_1,..., \tilde{v}_n\}$ be a basis of  $\mathcal{V}$. Then 
$\{\tilde{v}_1w,..., \tilde{v}_nw\}$ is a set of linearly independent functions, and, for each $i=1,...,n,$ there exists $(\tilde{m}_{ij})_{j=1,...n}$, such that
$$
- \frac{(\tilde{v}^iw)'}{w}\ \  =\ \ \sum_{j=1}^n \tilde{m}_{ij} \tilde{v}^j.
$$ 
It follows that $\mathbf{v}=(v_1,...,v_n):=(\tilde{v}_1 w, ..., \tilde{v}_nw)$ solves \eqref{Mstable} with
 $M=(m_{ij})$, $m_{ij}=-\tilde{m}_{ij}$. Moreover,  since  $\mathbf{v}$ solves \eqref{Mstable}, the $v^i$'s must be linear combination of functions  of the form  \eqref{str}, where the $\lambda_p$'s  are the eigenvalues of $M$. So, we also deduce that the eigenvalues of $M$ must have  real part strictly greater than $\lambda^*$, as $\tilde{v}_i\in L^{2}_w$ for all $i=1,...,n$, and that actually  $\tilde{v}_i\in W^{1,2}_w$. 
\hfill$\square$\\

In view of Proposition \ref{prop:stable}, we see that the $n$-dimensional $\mathcal{D}_w$-stable subspaces $\mathcal{V}$ of $L^2_w$ are in the form
\begin{multline}\label{structure}
\mathcal{V}\ \ =\ \ \mbox{Span}\  \Big\{w(\xi)^{-1}\xi^j e^{Re(\lambda_p)\xi}\cos(Im(\lambda_p)\xi),\ \  w(\xi)^{-1}\xi^j e^{Re(\lambda_p)\xi}\sin(Im(\lambda_p)\xi),\ 
\\   0\ \leq \ j\ \leq \ n_p-1, \  \ 1 \ \leq \  p \ \leq \ k\Big\},
\end{multline}
for some $k\geq 1$ and 
\be\label{constr}
(n_1, ..., n_k)\ \in\  \mathbb{N}^k \ \ \mbox{s.t.}\ \ n_1+...+n_k\ =\ n,\ \ \ \ \ \ \ 
 (\lambda_1,...,\lambda_k) \ \in\ \mathcal{C}^k, 
\ee
where
\begin{multline*}
 \mathcal{C}^k \ \ = \ \ \{z=(z_1,...,z_k)\in \mathbb{C}^k \ | \ z_i\neq z_j, \ \forall i\neq j; \  \ Re(z_j)>\lambda^*, \ \ \forall j=1,...,k; \\ \ \forall j\in\{1,...,k\}  \ \exists i\in\{1,...,k\} \mbox{ s.t. } \bar{z_i}=z_j \}.\end{multline*}
 Conversely,  all the subspaces $\mathcal{V}$ in the form \eqref{structure} above, with $(n_1, ..., n_k)$ and $(\lambda_1,...,\lambda_p)$ satisfying \eqref{constr},  are $n$-dimensional $\mathcal{D}_w$-stable subspaces of $L^2_w$. \\

Now, given an $n$-dimensional subspace $\mathcal{V}\subset L^2_w$,  denote
$$\bar{\mathcal{V}}\ \ :=\ \ \{ \mathbf{x}\in H_w \ | \ x_0=0, \ x_1\in \mathcal{V}\}.$$
\begin{definition}
We say that an $(n+1)$-dimensional subspace  $V\subset H_w $ is $\mathcal{A}^*$- stable if $V\subset D(\mathcal{A}^*)$ and $\mathcal{A}^*V\subset V$.
\end{definition}
{Noticing that $\mathcal{A}^* \mathbf{e}^0= \mathbf{0}$}, we immediately get the following corollary.

\begin{corollary}
An $(n+1)$-dimensional subspace $V\subset H_w $ is $\mathcal{A}^*$- stable if and only if
$$V\ \ =\ \ \mbox{Span} \ \{\mathbf{e}^0, \bar{\mathcal{V}}\},$$
with $\mathcal{V}$ being some $n$-dimensional $\mathcal{D}_w$-stable subspace of $L^2_w$.
\end{corollary}

\subsection{Exact finite-dimensional representation}\label{finite.sec}
Suppose that $\balpha,\bbeta\in V$, 
where $V$ is an $(n+1)$-dimensional $\mathcal{A}^*$-stable subspace of $H_w$. Let $\{\mathbf{e}^0, \mathbf{e}^1, ..., \mathbf{e}^n\}$ be an orthonormal basis of $(V, \langle\cdot,\cdot\rangle_w)$ and define
$$X^{k}_t=\langle \mathbf{e}^{k}, \mathbf{X}_{t}\rangle_w, \ \  \ \alpha^{k}=\langle \mathbf{e}^{k},\balpha\rangle_w, \ \  \ \beta^{k}=\langle \mathbf{e}^{k},\bbeta\rangle_w, \ \ \ \gamma^{k}=\langle \mathbf{e}^{k},\bgamma\rangle_w, \ \ \  k=0,...,n.$$
Then, since the projection of $\balpha, \bbeta$ onto $V^\perp$ is the null vector, we can write  the dynamics of $S_t= X^{0}_t$ as
\be\label{eq:app}
 dS_t\ \ = \ \ b\Big(S_t,\, \sum_{k=1}^n \alpha^{k}X^{k}_t,{ u_t}\Big)dt  +\sigma \Big(S_t,\, \sum_{k=1}^n \beta^{k}X^{k}_t,{ u_t}\Big)dW_t.
 \ee
 By $\mathcal{A}^*$-stability of $V$ we have the existence of a vector $\textbf{q}=(q_{k0})_{k=1,...,n}\in \mathbb{R}^n$ and of a matrix $\textbf{Q}=(q_{kh})_{h,k=1,...,n}\in \mathbb{R}^{n\times n}$ such that
\be
\mathcal{A}^*\mathbf{e}^{k}\ \ =\ \ q_{k0}\mathbf{e}^0+ \sum_{h=1}^nq_{kh}\mathbf{e}^{h}, \ \  k=1,..,n.
\ee

The dynamics of $S_t$ involves the processes $X^{k}(t)$, $k=1,...,n$, whose dynamics, plugging $\mathbf{e}^1, ..., \mathbf{e}^n$ in place of $\bvarphi$ in the definition of weak solution,   can be expressed in terms of themselves and of $S$ by means of the vector $\textbf{q}$ and of the matrix $\textbf{Q}$  as
 \be\label{eq:app2}
 dX^{k}_t\ \ =\ \ \Big(q_{k0}S_t+\sum_{h=1}^n q_{kh} X^{h}_t\Big)dt, \ \ \ \ k=1,...,n.
 \ee
The system of $n+1$ equations \eqref{eq:app} and \eqref{eq:app2} provides an $(n+1)$-dimensional Markovian representation of \eqref{SDDE} with initial datum $\mathbf{s}=(s_0,s_1(\cdot))$, the corresponding initial data being 
 \be\label{id}
 (x^0,x^{1},...,x^{n})\ \ =\  \ \left(\langle \mathbf{e}^{0},\mathbf{s}\rangle_w,\langle \mathbf{e}^{1},\mathbf{s}\rangle_w, ..., \langle \mathbf{e}^{n},\mathbf{s}\rangle_w\right).
 \ee
 {If also 
$$ 
\bgamma\ :=\ (0,\gamma_1)\in V,$$}
then the projection of $\bgamma$ onto $V^\perp$ is the null vector and  we can write 
$$Z_t\ \ =\ \ \sum_{k=0}^n \gamma^k X^k_t,$$
getting a representation of the process $Z$ in terms of the
$(n+1)$-dimensional (controlled)  Markov diffusion $(X^0,X^1, ...,
X^n)$.

\begin{remark}{From Proposition \ref{prop:stable}, it is clear that $\balpha,\bbeta$
and $\bgamma$ belong to an $(n+1)$-dimensional $\mathcal{A}^*$-stable
subspace of $H_w$ if and only if  the original coefficients $\tilde\balpha,\tilde\bbeta$
and $\tilde\bgamma$ belong to an $(n+1)$-dimensional $\mathcal{A}^*$-stable
subspace of $H$. Therefore, the property of having a
finite-dimensional Markovian representation does not depend on the
choice of the weight function. } 
\end{remark}
\subsection{Countable representation}
In general $\balpha,\bbeta, \bgamma$ fail to lie in an
$\mathcal{A}^*$-stable finite-dimensional subspace of $H_w$. However
one can consider an increasing sequence of  $\mathcal{A}^*$-stable
subspaces of $H_w$ and expand the problem along this sequence. In
order to construct such an increasing sequence, we  consider
specific subclasses of the general representation \eqref{structure}.  The simplest case to consider 
amounts to take $k=1$ in  \eqref{structure}, i.e. to consider a sequence of subspaces of $L^2_w$ in the form
\begin{equation}\label{structure2}
\mathcal{V}^n\ \ =\ \ \mbox{Span}\  \Big\{w(\xi)^{-1}\xi^j e^{\lambda\xi}, \ \ \ j=0,...,n-1\Big\}, \ \ \ \lambda>\lambda^*.
\end{equation}
To simplify  
the orthogonalization procedure of the basis of such subspaces,  { we restrict our analysis to the case
of exponential weights, i.e.  $w(\xi)= e^{p\xi}$, $p\in\mathbb{R}$}, { for which the
  construction of an orthonormal basis can be reduced to a known case,
  as we will see below}. {In this case clearly $\lambda^*=p/2$ and we
can choose, e.g.,  $\lambda>\max\{p,p/2\}$ to satisfy the constraint on $\lambda$  in \eqref{structure2}.
With this choice,  setting ${p}_0:=\lambda-p$ we have $p_0>0$, and the subspaces in \eqref{structure2} are rewritten as
\begin{equation}\label{structure3}
\mathcal{V}^n\ \ =\ \ \mbox{Span}\  \Big\{e^{{p}_0\xi} \xi^j, \ \ \ j=0,...,n-1\Big\}
\end{equation}}
or equivalently
\begin{equation}\label{structure4}
\mathcal{V}^n\ \ =\ \ \mbox{Span}\  \Big\{e^{{{p}_0}\xi} (2{p}_0\xi)^j, \ \ \ j=0,...,n-1\Big\}.
\end{equation}
The subspaces \eqref{structure4}, for running $n$, are the sequence of subspaces we shall consider.
An orthogonal basis with respect to the inner product {$\langle\cdot,\cdot\rangle_{L^2_w}$} for
the subspaces above can be constructed from the Laguerre polynomials as follows.
Define, for $k\geq 0$, the Laguerre poynomials  
$$
\tilde{P}_k (\xi)\ \ :=\ \ \sum_{i=0}^k \left(\begin{array}{c} k\\i\end{array}\right) \frac{(-1)^i}{i!}\xi^i, \ \ \ \ \xi\geq 0.
$$
Since we are working with $\mathbb{R}^-$ instead of $\mathbb{R}^+$ the
sign of the argument is inverted so we consider Laguerre's polynomials on $\mathbb{R}^-$ defined as
$$
P_k(\xi)\ \ =\ \ \tilde{P}_k(-\xi), \ \  \ \ \xi\leq 0.
$$
From the definition of $P_k$'s, we have
\be\label{derp}
P_k(0)\ \ =\ \ 1, \ \  \ k\geq 0,
\ee
and moreover, using an induction argument,   we get
\be\label{derp1}
P_k'\  \ =\ \ \sum_{i=0}^{k-1} P_i, \ \ \  k\geq 1.
\ee
As is well known, the Laguerre functions 
$$
l_{k,{p}_0}(\xi) \ \ = \ \ \sqrt{2{p}_0}\ P_k(2{p}_0\xi)\ e^{{p}_0\xi}
$$
are an orthonormal basis for $L^2$. So the functions
\begin{equation}\label{Lag}
L_{k,{p}_0,p}(\xi)\ \ =\ \ e^{-\frac{p}{2}\xi}\  l_{k,{p}_0}(\xi),
 \ \ \ \ \ \ k=0,...,n-1, 
\end{equation}
are an orthonormal basis with respect to the inner product { $\langle\cdot,\cdot\rangle_{L^2_w}$} for $\mathcal{V}^n$ defined in \eqref{structure4} and the sequence of functions
$$
(L_{k,{p}_0,p})_{k\geq 0}
$$
is an orthonormal basis for $L^2_w$.

Consider the system of vectors  $(\mathbf{e}_k)_{k\geq 0}$ in $H_w$ where
$$\mathbf{e}^0 \ =\  (1,0); \ \ \ \ \mathbf{e}^k\ =\ (0,L_{k-1,{p}_0,p}), \ \ k\geq 1.$$
Then, from the argument above, this system is  an orthonormal basis  in $H_w$. 

{
Using \eqref{derp}--\eqref{Lag}, we have for $k\geq 1$ (with the convention that $\sum_{i=1}^0 = 0$)
\begin{align*}
L'_{k-1,{p}_0,p} (\xi)&= \sqrt{2{p}_0}\frac{d}{d\xi} (e^{({p}_0-p/2)\xi} P_{k-1}(2{p}_0\xi))\\
& =  ({p}_0-p/2)
L_{k-1,{p}_0,p} (\xi)+ (2{p}_0)^{3\over 2} e^{({p}_0-p/2)\xi}\sum_{i=0}^{k-2}
P_i(2{p}_0\xi)\\
& = ({p}_0-p/2)
L_{k-1,{p}_0,p} (\xi)+ 2{p}_0 \sum_{i=0}^{k-2}L_{i,{p}_0,p} (\xi)
\end{align*}}
So, by \eqref{{A}^*} 
\be\label{A*e}
\mathcal{A}^*\mathbf{e}^0\ \ = \ \ \mathbf{0};  \ \ \ \ \ \  \ \
\mathcal{A}^*\mathbf{e}^k\ \ =\ \ {\mathbf{e}^0 \sqrt{2{p}_0}-2{p}_0\sum_{i=1}^{k-1} \mathbf{e}^i-{  ({p}_0-p/2)}\,\mathbf{e}^k,} \ \ \ k\geq 1.
\ee
\begin{remark}
From \eqref{A*e} we see that, 
setting 
$$
V^n\ :=\ \mbox{\mbox{Span}}\, \left\{\mathbf{e}^0, ... , \mathbf{e}^n\right\}, \ \ \ \ \ \ n\geq 0,
$$
we have the $\mathcal{A}^*$-stability of $V^n$ for each $n\geq 0$. 
\end{remark}
\medskip
Setting
$$X^{k}(t)\  =\ \langle \mathbf{e}^{k}, \mathbf{X}_{t}\rangle_w, \ \ \ \ \  \ \alpha^{k}\  = \ \langle \mathbf{e}^{k},\balpha\rangle_w, \ \  \ \ \ \  \beta^{k} \ = \ \langle \mathbf{e}^{k},\bbeta\rangle_w, \ \ \  k\geq 0,$$ we have the Fourier series expansions in $H_w$
 $$\mathbf{X}_{t}\   = \ \sum_{k=0}^\infty X^k(t)\mathbf{e}^k, \ \ \ \ \ 
 \balpha\  =\ \sum_{k=0}^\infty\alpha^k\mathbf{e}^k, \ \ \ \ \
 \bbeta\  = \ \sum_{k=0}^\infty\beta^k\mathbf{e}^k.$$
 Then we can rewrite the dynamics of $S_t=X^0_t$ as
 \be\label{count}
 dS_t\ \ =\ \ b\,\Big(S_t,\, \sum_{k=1}^\infty \alpha^kX^k_t,\,u_{t}\Big)dt  \ +\ \sigma \,\Big(S_t,\, \sum_{k=1}^\infty \beta^kX^k_t,\,u_{t}\Big)dW_t.
 \ee
Using the definition of weak solution Definition \ref{def:mild-weak}-(ii) and considering also \eqref{A*e} we have
 \be\label{count2}
 dX^k(t)\ \ =\ \ {{\Big(\sqrt{2{p}_0}S_t- 2{p}_0\sum_{i=1}^{k-1} X^i_t-  ({p}_0-p/2) X^k_t\Big)dt,}} \ \ \ \ \ k\geq 1.
\ee
Setting the  initial data
 \be\label{idinf}
\big(s_0,(x^{k})_{k\geq 1}\big)\ \ =\ \ \big(\langle \mathbf{e}^{0},\mathbf{s}\rangle_w,(\langle\mathbf{e}^k,\mathbf{s}\rangle_w)_{k\geq 1}\big),
 \ee
equations \eqref{count}--\eqref{count2} provide a  countable
Markovian representation  of our original system \eqref{SDDE}.
Moreover, setting also
$$\gamma^k \   \ =\ \ \ \langle \mathbf{e}^{k}, \bgamma\rangle_w,$$
we have the  Fourier series expansion for $\bgamma$
$$
\bgamma\ \ =\ \  \sum_{k=0}^\infty \gamma^k\mathbf{e}^k,
$$
so we also have the representation of the process \eqref{function1} as
\be\label{fun1}
Z_t\ \ =\ \ \sum_{k=0}^\infty \gamma^kX^k(t) \ \ =\ \ \gamma^0S_t+\sum_{k=1}^\infty \gamma^k X^k(t).
\ee
\subsection{Approximate finite-dimensional  representation}
When  $\balpha,\bbeta$ belong to some finite dimensional subspace
$V^n$,  equations \eqref{count}-\eqref{count2}-\eqref{idinf} provide a finite-dimensional representation of \eqref{SDDE} in the spirit of the previous subsection. In this case the dynamics of $S$ requires only the knowledge of $(X^k)_{k=0,...,n}$ and the dynamics of these variables is given also in terms of themselves. Finally, when also $\bgamma$ belongs to $V^n$, then  
\eqref{function1} can be written in terms of the finite-dimensional
Markov process $(X^k)_{k=0,...,n}$ and we fall into an exact
finite-dimensional representation of the problem. When some of the
above conditions fail to be true (i.e. there is no finite dimensional subspace $V^n$ such that $\balpha,\bbeta,\bgamma\in V^n$), then we need to truncate the Fourier series for $\balpha, \bbeta, \bgamma$ and work with an approximate finite-dimensional representation of the problem. In this case,
setting for $n\geq	0$
$$
\balpha^n=\sum_{k=0}^n \alpha^k\mathbf{e}^k, \ \ \ \ \ \bbeta^n=\sum_{k=0}^n \beta^k\mathbf{e}^k, \ \  \ \ \ \ \bgamma^n=\sum_{k=0}^n \gamma^k\mathbf{e}^k,
$$
 we have the following estimate for the error.

\begin{proposition}\label{err.prop}
For each $n\geq 0$  let $(S^n,X^{n,k})_{k=0,...,n}$ be the finite dimensional Markov diffusion solving
 \be\label{count0}
 dS^n_t\ \ =\ \ b\,\Big(S^n_t,\, \sum_{k=1}^n \alpha^kX^{n,k}_t,\,u_t\Big)dt  \ +\ \sigma \,\Big(S^n_t,\, \sum_{k=1}^n \beta^kX^{n,k}_t,\,u_t\Big)dW_t.
 \ee
 \be\label{count20}
 dX^{n,k}_t\ \ =\ \ { \Big(\sqrt{2{p}_0}S^n_t- 2{p}_0\sum_{i=1}^{k-1} X^{n,i}_t-  ({p}_0-p/2)X^{n,k}_t\Big)dt,} \ \ \ \ \ k=1, ..., n,
\ee
with initial data
 \be\label{idinf1}
\big(s_0,(x^{n,k})_{k=1,...,n}\big)\ \ =\ \ \big(\langle \mathbf{e}^{0},\mathbf{s}\rangle_w,(\langle\mathbf{e}^k,\mathbf{s}\rangle_w)_{k=1,...,n}\big).
 \ee
{Then for every $T>0$, there exists
$C=C_{T,\|\mathbf{s}\|_w,\|\salpha\|_w,\|\sbeta\|_w,\|\sgamma\|_w}<\infty$
such that, uniformly on $u\in\mathcal{U}$, 
$$
\mathbb{E}\Bigg[\sup_{0\leq t\leq T}\Big|\,Z_t-\sum_{k=0}^n\gamma^k X^{n,k}_t\,\Big|^2\Bigg]\ \ \leq\ \  C(\|\balpha-\balpha^n\|_w^2+\|\bbeta-\bbeta^n\|_w^2+\|\bgamma-\bgamma^n\|_w^2).
$$}
\end{proposition}
\textbf{Proof.} {From Proposition \ref{prop:growthX} it follows that
$S_t$, $S^n_t$, $X^k_t$ and $X^{n,k}_t$ are square integrable for all
$t\geq 0$, all $n\geq 1 $ and all $k\in\{1,...,n\}$. 
Using standard tools such as Doob's inequality and It\^o isometry, one can show firstly that
\begin{align*}
&\mathbb E[\sup_{0\leq t\leq T}(S_t - S^n_t)^2]  \leq 2T\int_0^T \mathbb E[(b(S_r,
\sum_{k=1}^\infty \alpha^k X^k_r, u_r) - b(S^n_r,
\sum_{k=1}^n \alpha^k X^{n,k}_r, u_r))^2] \\ &\qquad \qquad+
8\int_0^T \mathbb E[(\sigma(S_r,
\sum_{k=1}^\infty \beta^k X^k_r, u_r) - \sigma(S^n_r,
\sum_{k=1}^n \beta^k X^{n,k}_r, u_r))^2] \\
& \leq 48 (T+1) C_1^2 \int_0^T \mathbb E\left\{(S_r-S^n_r)^2 +
\left(\sum_{k=n+1}^\infty \alpha^k X^k_r\right)^2 +
\left(\sum_{k=n+1}^\infty \beta^k X^k_r\right)^2+\left(\sum_{k=1}^n
  \alpha^k (X^k_r-X^{n,k}_r)\right)^2\right\}dr<\infty
\end{align*}
and similarly,
$$
\mathbb E[\sup_{0\leq t \leq T}(X^k_t - X^{n,k}_t)^2]<\infty,\quad k=1,\dots,n.
$$
Then, let us introduce the quantity
$$
{ M_T} \ := \ \mathbb E[\sup_{0\leq t \leq T}(S_t - S^n_t)^2] + \sum_{k=1}^n \mathbb E[\sup_{0\leq t \leq T} (X^k_t -
X^{n,k}_t)^2]\ <\ \infty. 
$$
From the above estimates, we then get
\begin{align*}
\mathbb E[\sup_{0\leq t\leq T}(S_t - S^n_t)^2]\ \leq \ 48 (T+1) C_1^2 \int_0^T \{(1+\|\balpha^n\|^2)M_r + (\|\balpha - \balpha^n\|^2_w + \|\bbeta - \bbeta^n\|^2_w)\mathbb E[\|\mathbf X_r\|^2_w]\}dr 
\end{align*}
and also for $k=1,\dots,n$,
\begin{align*}
\mathbb E[\sup_{0\leq t\leq T}(X^k_t -
X^{n,k}_t)^2] \ \leq \ (k+1)(4p^2+2p)T \int_0^T M_r dr,
\end{align*}
so that for some constant $C$ depending on $p$, $n$, $T$ and
$\alpha$, 
$$
M_T\  \leq \ C \int_0^T \{M_r + (\|\balpha - \balpha^n\|^2_w + \|\bbeta - \bbeta^n\|^2_w)\mathbb E[\|\mathbf X_r\|^2_w]\}dr. 
$$
From Gronwall's inequality and Proposition \ref{prop:growthX} it follows that there exists another constant,
also denoted by $C$, depending on $p$, $n$, $T$,
$\alpha$, $\beta$ and the initial condition, such that 
$$
M_T \ \leq \ C(\|\balpha - \balpha^n\|^2_w + \|\bbeta - \bbeta^n\|^2_w). 
$$
Finally, 
\begin{align*}
\mathbb{E}\Bigg[\sup_{0\leq t\leq T}\Big|\,Z_t-\sum_{k=0}^n\gamma^k X^{n,k}_t\,\Big|^2\Bigg]
&\ \leq \ 2
\mathbb{E}\Bigg[\sup_{0\leq t \leq T}\Big|\sum_{k=n+1}^\infty\gamma^k
X^{k}_t\Big|^2\Bigg]+\mathbb{E}\Bigg[\sup_{0\leq t \leq T}\Big|\sum_{k=0}^n\gamma^k (X^k_t -
X^{n,k}_t)\Big|^2\Bigg]\\
&\ \leq \ 2\|\bgamma-\bgamma^n\|_w^2 \mathbb E[\sup_{0\leq t \leq T}\|\mathbf
X_t\|^2_w] + \|\bgamma^n\|^2_w M_T.  
\end{align*}
Combining this with the bound on $M_T$ obtained above and Proposition \ref{prop:growthX}, the proof is
complete. }
\hfill$\square$\\

This proposition shows that the error of approximating the process $Z$
with a linear combination of components of a finite-dimensional
Markovian diffusion $\sum_{k=0}^n \gamma^k X^{n,k}_t$ depends on the
error of approximating the coefficients  $\balpha$, $
\bbeta$ and $\bgamma$ with the corresponding truncated Fourier-Laguerre
series. The actual convergence rate as $n\rightarrow  \infty$ will depend on
the regularity of the functions $\balpha$, $\bbeta$ and
$\bgamma$. For example, from \cite[Lemma A.4]{BTW10} it follows that
if these functions are constant in the neighborhood of zero, have
compact support and finite variation (this is the case e.g., for
uniformly weighted moving averages) and $w\equiv 1$ then
$\|\balpha-\balpha^n\|_w^2+\|\bbeta-\bbeta^n\|_w^2+\|\bgamma-\bgamma^n\|_w^2
\leq C n^{-3/2}$ for some constant $C$ and $n$ sufficiently large. For
$C^\infty$ functions, on the other hand, the convergence rates are
faster than polynomial. 
\section{Application to optimal control and  stopping }\label{app.sec}
{In this  last section we show how the results of the previous can be implemented to treat optimal control or optimal stopping problems. { Within this section it is assumed that $S$ solves \eqref{SDDE1}-\eqref{IC} and $Z$ is the process defined in \eqref{function1}}.}
\subsection{Optimal control problems}
Let $T>0$. Given $f:[0,T]\times \mathbb{R}\rightarrow\mathbb{R},\ \phi:\mathbb{R}\rightarrow\mathbb{R}$, we may  consider following optimal control  problem 
$$V(\mathbf{s})\ \ =\ \ \inf_{u\in\mathcal{U}} \ J(\mathbf{s};u),$$
where
$$
J(\mathbf{s};u) \ :=\ \mathbb{E}\left[\int_0^Tf(t,Z_t,{ u_t})dt+\phi(Z_T)\right].$$
This problem cannot be solved by Dynamic Programming due to the lack of markovianity. However, given $n\geq 0$, we can consider the problem in $\mathbb{R}^{n+1}$
$$V^n(\mathbf{x}^n)\ \ =\ \ \inf_{u\in\mathcal{U}}\ J^n(\mathbf{x}^n;u),\ \ \ \ \mathbf{x}^n=(x^n_0,x^n_1,...,x^n_n)\in\mathbb{R}^{n+1},$$
where
$$J^n(\mathbf{x}^n;u)\ :=\ \mathbb{E}\left[\int_0^Tf(t,Z^{n}_t,{ u_t})dt+\phi(Z^{n}_T)\right],$$
and
\begin{itemize}
\item[(i)] the ``output" process $Z^{n}_T$ is 
\begin{align}
Z^{n}_t\ \ =\ \ \sum_{k=0}^n\gamma_kX^{n,k}_t;\label{system2z}
\end{align}
\item[(ii)]
the state equation for the $(n+1)$-dimensional process $(X^{n,k})_{k=0,...,n}$ is 
\be\label{system2}
\begin{cases}
 dX^{n,0}_t\ =\  b\left(X^{n,0}_t,\, \sum_{k=1}^n \alpha^kX^{n,k}_t,\,u_{t}\right)dt +\sigma \left(X^{n,0}_t,\, \sum_{k=1}^n\beta^kX^{n,k}_t,\,u_{t}\right)dW_t,\\\\
 dX^{n,k}_t\ =\ {\left(\sqrt{2{p}_0}X^{n,0}_t- 2{p}_0\sum_{i=1}^{k-1}X^{n,i}_t
 -  ({p}_0-p/2)X^{n,k}_t\right)dt}, \ \ \ \ k=1,...,n,
 \ec
 \ee
with initial data
\be\label{ind}
X^{n,k}_0\ \ =\ \ x^n_k, \ \ \ \ k=0,...,n.
\ee
\end{itemize}
This finite-dimensional problem can be solved via the corresponding
Hamilton-Jacobi-Bellman equation \cite{FS06}. The following
proposition provides an error estimate for the value function. 
\begin{proposition}\label{prop:value}
Suppose that $f(t,\cdot,u)$ is Lipschitz continuous uniformly in $t\in[0,T]$ { and $u\in U$}, and that $\phi$ is Lipschitz continuous.
Set 
$$ \mathbf{x}^n(\mathbf{s})\ \ =\ \ (\langle \mathbf{s},\mathbf{e}^k\rangle_w)_{k=0,...,n}, \ \ \ n\geq 0.$$ 
Then there exists $K=K_{T,\|\mathbf{\mathbf{s}}\|_w,\|\salpha\|_w,\|\sbeta\|_w,\|\sgamma\|_w}$ such that 
$$|V(\mathbf{s})-V^n(\mathbf{x}^n(\mathbf{s}))|^2\ \ \leq\ \  K(\|\balpha-\balpha^n\|^2+\|\bbeta-\bbeta^n\|^2+\|\bgamma-\bgamma^n\|^2), \ \ \ \ \ \forall n\geq 0.$$
\end{proposition}
\textbf{Proof.}
{We shall use the notation $Z^{u,\mathbf s}$ and $Z^{n,u,\mathbf x^n(\mathbf{s})}$
to make the dependence on the initial condition and the control
explicit. The common Lipschitz constant of $f$ and $\phi$ shall be
denoted by $K_0$.
We have
{
\begin{align*}
 |V(\mathbf{s})-V^n(\mathbf{x}^n(\mathbf{s}))| &\leq  \sup_{u\in \mathcal U}
|J(u) - \inf_{u\in \mathcal U} J^n(u)| \\
& \leq  \sup_{u\in \mathcal U} \mathbb E\left[\int_0^T |f(t,Z^{u,\mathbf
    s}_t,{ u_t}) - f(t,Z^{n,u,\mathbf x^n(\mathbf
    s)}_t,{ u_t})|dt + |\phi(Z^{u,\mathbf s}_T) -\phi(Z^{n,u,\mathbf x^n(\mathbf
    s)}_T)| \right]\\
&\leq  K_0 \sup_{u\in \mathcal U} \mathbb E\left[\int_0^T |Z^{u,\mathbf
    s}_t - Z^{n,u,\mathbf x^n(\mathbf
    s)}_t|dt + |Z^{u,\mathbf s}_T -Z^{n,u,\mathbf x^n(\mathbf
    s)}_T|\right]\\
&\leq K_0\sup_{u\in\mathcal U} \left\{\int_0^T \mathbb E[|Z^{u,\mathbf
    s}_t - Z^{n,u,\mathbf x^n(\mathbf
    s)}_t|^2]^{\frac{1}{2}}dt + \mathbb E[|Z^{u,\mathbf
    s}_T - Z^{n,u,\mathbf x^n(\mathbf
    s)}_T|^2]^{\frac{1}{2}}\right\}\\
&\leq K_0 C^{\frac{1}{2}}_T (T+1) (\|\balpha - \balpha^n\|^2_w+\|\bbeta - \bbeta^n\|^2_w+\|\bgamma - \bgamma^n\|^2_w)^{\frac{1}{2}},
\end{align*}
where the last inequality follows from Proposition \ref{err.prop}, and $C_T$ is the bound on the constant $C$ of that proposition over $t\in[0,T]$.} 
}
\hfill$\square$\\\\
The result above can be applied, for instance, to the problem investigated in \cite{GM06,GMS09} or to generalizations of the examples shown in \cite{EOS00}.
\subsection{Optimal stopping problems}\label{Sub:optstop}
Let $T>0$ and consider \eqref{SDDE} when $b,\sigma$ do not depend on $u$ (so the diffusion is actually uncontrolled). Letting  $\mathcal{T}$ be the set of all stopping times with respect to the filtration $(\mathcal{F}_t)_{t\in[0,T]}$ and taking values in the interval $[0,T]$ and given $\phi:[0,T]\times \mathbb{R}\rightarrow\mathbb{R}$, we may  consider following optimal stopping problem 
$$V(\mathbf{s})\ \ =\ \ \inf_{\tau\in\mathcal{T}}\mathbb{E}\left[\phi(\tau,Z_\tau)\right].$$
Also in this case the problem cannot be solved by Dynamic
Programming due to the lack of markovianity. However, to approximate
its solution, given $n\geq 0$, we can consider the problem in $\mathbb{R}^{n+1}$
$$V^n(\mathbf{x}^n)\ \ =\ \ \inf_{\tau\in\mathcal{T}}\ \mathbb{E}\left[\phi(\tau,Z^{n}_\tau)\right],\ \ \ \ \mathbf{x}^n=(x^n_0,x^n_1,...,x^n_n)\in\mathbb{R}^{n+1},$$
where the ``output'' process $Z^{n}_t$ and the state equation for the
$(n+1)$-dimensional process $(X^{n,k})_{k=0,...,n}$ are given,
respectively, by \eqref{system2z} and \eqref{system2}.
The following proposition provides an error estimate. 
\begin{proposition}
Suppose that $\phi(t,\cdot)$ is Lipschitz continuous uniformly in $t\in[0,T]$.
Set 
$$ \mathbf{x}^n(\mathbf{s})\ \ =\ \ (\langle \mathbf{s},\mathbf{e}^k\rangle_w)_{k=0,...,n}, \ \ \ n\geq 0.$$ 
Then there exists $K=K_{T,\|\mathbf{\mathbf{s}}\|_w,\|\salpha\|_w,\|\sbeta\|_w,\|\sgamma\|_w}$ such that 
$$|V(\mathbf{s})-V^n(\mathbf{x}^n(\mathbf{s}))|^2\ \ \leq\ \  K(\|\balpha-\balpha^n\|^2+\|\bbeta-\bbeta^n\|^2+\|\bgamma-\bgamma^n\|^2), \ \ \ \ \ \forall n\geq 0.$$
\end{proposition}
\textbf{Proof.} Similar to the proof of Proposition \ref{prop:value}.\hfill$\square$\\

This result can be applied, for example, to the problem of pricing
American options written on the moving average of the asset price,
which was studied in \cite{BTW10} by approximating the dynamics by a
finite-dimensional  Markovian one by means of Laguerre polynomials,
but without passing through the infinite-dimensional representation of
the system. Let us briefly recall the problem. 
Let $T>0$ and let $S$ be the price of a stock index and consider the financial problem of pricing an American option whose payoff at the exercise time $t\in [0,T]$  depends on a past average of the stock, i.e. 
$$\phi\left(\frac{1}{\delta}\int_{t-\delta}^tS_{\xi}d\xi\right), \ \ \delta>0.$$
Suppose that the price $S$ is a Markov diffusion solving the SDE
 \be\label{SDE}
dS_t\ \ =\ \ b\left(S_t\right)dt+ \sigma\left(S_t\right)dW_t, \ \ \ \ \ S_0\ =\ s_0\ >\ 0,
\ee
and set
$$Z_t\ \ =\ \ \frac{1}{\delta}\int_{t-\delta}^tS(r)dr\ \ =\ \ \frac{1}{\delta}\int_{-\delta}^0S_{t+\xi}d\xi.$$
Letting $\mathcal{T}$ be the set of all stopping times with respect to the filtration $(\mathcal{F}_t)_{t\in[0,T]}$ taking values in the interval $[0,T]$, 
the value of the option at time $0$ is
$$
V(s_0,s_1)\ \ =\ \ \sup_{\tau\in\mathcal{T}}\mathbb{E}\left[\phi(Z(\tau))\right],
$$
where $s_1(\xi)$  is the value of the stock at the past time $\xi\in[-\delta,0]$.
This problem is intrinsically infinite-dimensional; it falls into our setting as a special case by taking
 $$\balpha,\bbeta\ =\ (1,0)\in H_w,\ \ \ \ \ \bgamma \ =\ \left(0,\
   \frac{1}{\delta}\mathbf{1}_{[-\delta,0]}\right)\in H_w.$$
For more details about this problem we refer to \cite{BTW10}.\\\\

{\textbf{Acknowledgements.} The authors thank Mauro Rosestolato for very fruitful discussions and suggestions.}


\end{document}